\documentclass[12pt]{article} %

\usepackage[style=authoryear-comp,firstinits=true,citestyle=authoryear,natbib=true,useprefix=true]{biblatex}

\DeclareNameAlias{sortname}{last-first}

\renewbibmacro{in:}{}
\usepackage{tikz}
\usepackage{pgfplots}
\pgfplotsset{compat=1.18}
\usepackage{enumitem}

\usepackage{array,epsfig,fancyheadings,rotating}

\usepackage[]{hyperref}  %
\usepackage{sectsty, secdot}
\sectionfont{\fontsize{12}{14pt plus.8pt minus .6pt}\selectfont}
\renewcommand{\theequation}{\thesection\arabic{equation}}
\subsectionfont{\fontsize{12}{14pt plus.8pt minus .6pt}\selectfont}

\usepackage{amsmath}
\usepackage{amssymb}
\usepackage{amsfonts}
\usepackage{multirow}
\usepackage{amsthm}

\newtheorem{theorem}{Theorem}
\newtheorem{lemma}{Lemma}
\newtheorem{corollary}{Corollary}
\newtheorem{proposition}{Proposition}
\theoremstyle{definition}
\newtheorem{definition}{Definition}
\newtheorem*{assumptions*}{Assumptions}
\newtheorem{example}{Example}
\newtheorem{remark}{Remark}
\usepackage{bm}

\DeclareMathOperator{\cv}{CV}
\DeclareMathOperator{\scv}{SCV}
\DeclareMathOperator{\ascv}{ASCV}
\DeclareMathOperator{\optimal}{Optimal}
\DeclareMathOperator{\var}{Var}

\DeclareMathOperator{\ive}{\overline{IV}}
\DeclareMathOperator{\ibsh}{\widehat{IBS}}
\DeclareMathOperator{\argmin}{argmin}
\DeclareMathOperator{\N}{N}

\DeclareMathOperator{\sure}{SURE}
\DeclareMathOperator{\mise}{MISE}
\DeclareMathOperator{\ise}{ISE}
\DeclareMathOperator{\vecspan}{span}
\newcommand{\bone}{\bm{1}}

\usetikzlibrary{external}
\usepackage{caption}
\DeclareCaptionLabelFormat{boldcap}{\bf #1 #2}
\DeclareCaptionLabelFormat{boldcap}{\bf #1 #2}
\begin{document}

\renewcommand{\baselinestretch}{2}

\markright{ \hbox{\footnotesize\rm 
}\hfill\\[-13pt]
\hbox{\footnotesize\rm
}\hfill }

\markboth{\hfill{\footnotesize\rm FIRSTNAME1 LASTNAME1 AND FIRSTNAME2 LASTNAME2} \hfill}
{\hfill {\footnotesize\rm FILL IN A SHORT RUNNING TITLE} \hfill}

\renewcommand{\thefootnote}{}
$\ $\par

\fontsize{12}{14pt plus.8pt minus .6pt}\selectfont \vspace{0.8pc}
\centerline{\large\bf STABILIZED CROSS-VALIDATION OF}
\vspace{2pt} 
\centerline{\large\bf SMOOTHNESS IN DENSITY DECONVOLUTION}
\vspace{.4cm} 
\centerline{David Kent}
\vspace{.4cm} 
\centerline{\it Cornell University}
 \vspace{.55cm} \fontsize{9}{11.5pt plus.8pt minus.6pt}\selectfont

\begin{quotation}
\noindent {\it Abstract:} \\

We consider density estimation under measurement error with the Smoothness-Penalized Deconvolution (SPeD) estimator.
The estimator has a tuning parameter regulating the smoothness of the estimate, and proper choice of this parameter is critical for forming good estimates.
We derive the cross-validation choice of tuning parameter for the SPeD estimator, but it performs very poorly.
We introduce a stabilized cross-validation (SCV) criterion which unbiasedly estimates the mean integrated squared error (MISE) for a smaller sample size, and use asymptotic arguments to obtain an appropriate tuning parameter from this stabilized criterion.
We show that the SCV is a strongly consistent estimator of the MISE, and that it is the minimum variance unbiased estimator of the MISE.
In a simulation study, we show that the SCV approach outperforms the previously recommended choice of tuning parameter in nearly all settings, and in a majority of the settings, SPeD with the SCV outperforms the classic deconvoluting kernel estimator with its recommended choice of tuning parameter.

\vspace{9pt}
\noindent {\it Key words and phrases:}
Penalty parameter, deconvolution, measurement error, cross-validation, density estimation
\par
\end{quotation}\par

\def\thefigure{\arabic{figure}}
\def\thetable{\arabic{table}}

\renewcommand{\theequation}{\thesection.\arabic{equation}}

\fontsize{12}{14pt plus.8pt minus .6pt}\selectfont

\tableofcontents %

\section{Introduction}
\label{sec:intro}

Suppose \(X\) has unknown pdf \(f\), \(E\) has known pdf \(g\), and \(X\) and \(E\) are independent.
Then \(Y = X + E\) has pdf \(h(y) = g*f(y) = \int f(t)g(y-t)\,dt\).
Suppose we observe a sample \(Y_1,\dots,Y_n \overset{iid}{\sim} h\), and we wish to estimate the density \(f\) of \(X\).
This is the deconvolution problem.
The recent Smoothness-Penalized Deconvolution (SPeD) estimator (\cite{kentSmoothnessPenalizedDeconvolutionSPeD2022}; \cite{yang_bias-corrected_2021}; \cite{yang_density_2020}) depends critically on a penalty parameter \(\alpha\) that regulates the smoothness of the estimate.
If \(\alpha\) is too large, the estimate is too smooth, and interesting features of \(f\) may be absent; if \(\alpha\) is too small, the estimate is too wiggly, and contains many features not present in \(f\).

In \cite{yang_density_2020}, the authors choose \(\alpha\) by a \emph{SURE} criterion, similar to Stein's unbiased risk estimation.
Specifically, they consider a new hypothetical sample \(Y_1^0, \dots, Y_n^0\) identically distributed as and independent from the original sample, and denote by \(\hat h\) and \(\hat h^0\) density estimates of \(h\) from the original and hypothetical samples, respectively.
Letting \(\hat f_{\alpha}\) denote the estimator fit to the original sample, they choose \(\alpha\) by minimizing an estimate of the SURE criterion \(\sure(\alpha) = \mathbb{E}\|\hat h^0 - g*\hat f_{\alpha}\|^2\), where the expectation is over both the original and hypothetical samples, which determine \(\hat f_{\alpha}\) and \(\hat h^0\), respectively.

The authors note that this SURE approach occasionally chooses an \(\alpha\) that is far too small, leading to a very rough estimate.
In our simulations in Section~\ref{sec:simulations}--see, for example, Figure~\ref{fig:muchworse}--we will see that these too-small choices of \(\alpha\) can lead to catastrophically poor estimates, which poses difficulties in practice.
This is not a huge surprise in retrospect: the deconvolution problem is ill-posed \citep[Section 2]{kentSmoothnessPenalizedDeconvolutionSPeD2022} in the sense that for two estimates \(\hat f\) and \(\hat f'\), the values \(g*\hat f\) and \(g*\hat f'\) can be very close even when \(\hat f\) and \(\hat f'\) are very different.
But since \(\sure(\alpha)\) measures the error in \(g*\hat f_{\alpha}\), this criterion may be insensitive to choice of \(\alpha\) even when the error in \(\hat f_{\alpha}\) is quite sensitive.
This observation aligns with the observation in \cite{yang_density_2020} that the problem occurs when there is a ``wide range of nearly-optimal [\(\alpha\)] values,'' indicating that the problem occurs when the error criterion is insensitive to the \(\alpha\) values.

Cross-validation has been suggested for choosing the bandwidth in the deconvolution kernel density estimator (\cite{stefanski_deconvoluting_1990}; \cite{hesse_data-driven_1999}; \cite{youndjeOptimalBandwidthSelection2008}), though the intuition behind it is less transparent in deconvolution than in classical density estimation.
In Section~\ref{sec:criterion:cv} we translate the intuitive justification for cross-validation into the deconvolution problem, and derive the cross-validation criterion for the SPeD estimator.
We find that it suffers from extremely poor performance, regularly choosing \(\alpha\) too small.
In fact, many realizations of \(\cv(\alpha)\) are monotonically increasing on any reasonable interval we might optimize over (see an example plotted in Figure~\ref{fig:realizations}); in these cases, the minimizer of \(\cv(\alpha)\) is determined by the largely-arbitrary lower bound of the optimization interval.

In classical density estimation, the cross-validation criterion consists of a non-random term that closely approximates the integrated variance, and a random term that estimates the integrated squared bias (IBS); the same is true in deconvolution.
In \cite{hall_smoothed_1992}, explaining the large variability of cross-validated bandwidths in the classical case, the authors note that ``estimates of integrated squared bias tend to suffer from high variance when bandwidths of the order \(h^{-1/5}\) are considered''---but of course this is the order of the optimal bandwidth, so it is exactly the order of the bandwidths we should be searching in cross-validation.
We see a similar phenomenon in practice: on the scale of the penalty parameter minimizing the MISE, the component of \(\cv(\alpha)\) estimating the bias has high variance.

This insight makes us wish we were considering \(\alpha\) of a larger order.
In Section~\ref{sec:criterion:scv} we do exactly that, introducing a procedure we call \emph{stabilized cross-validation} (SCV).
In SCV, we adjust the CV criterion to estimate the error for a sample of size \(m \ll n\); this adjusted criterion consists of a non-random term closely approximating the integrated variance (for a sample of size \(m\)), plus \emph{the same random term} (to a very close approximation) that occurs in the CV criterion approximating the IBS.
For a given \(\alpha\), the variance of the IBS is essentially unchanged, but when we deal with the SCV criterion we are considering \(\alpha\) of a larger order than with the CV criterion due to the increased integrated variance term.
For \(m \approx n^{1/2}\), the SCV criterion is much better behaved, and its minimum has much lower variance.
After minimizing the SCV criterion to estimate \(\alpha_m'\) suitable for a sample of size \(m\), we use an idea from \cite{hallUsingBootstrapEstimate1990} to parlay this into an estimate \(\hat \alpha_n\) suitable for a sample of size \(n\).
Namely, we use \emph{a priori} knowledge of the asymptotic rate of the optimal \(\alpha_n\), say \(\alpha_n \sim \kappa b_n\) for a known sequence \(b_n\) and unknown constant \(\kappa\), and form \(\hat \alpha_n = (b_n/b_m) \alpha_m'\).

\subsection{The Estimator}

The estimator we address is a slight modification of the SPeD addressed in \cite{kentSmoothnessPenalizedDeconvolutionSPeD2022}, which will simplify its analysis.
In their definition of the estimator, they begin with a preliminary density estimate \(h_n\) of \(h\); here, we replace \(h_n\) by the empirical distribution, which amounts to replacing \(\int [g*u]h_n\) by \(\frac1n \sum_{j=1}^n g*u(Y_i)\).
For a similar approach in a different problem, see \cite{wager_geometric_2014}.
Here we present three equivalent representations of this estimator (cf. Theorem 2 of \cite{kentSmoothnessPenalizedDeconvolutionSPeD2022}).
We will denote the Fourier transform by an overset twiddle, so that \(\tilde f(t) = \mathbb{E}[e^{-itX}] = \int e^{-itx}f(x)\,dx\). Let \(P_n = \frac1n\sum_{j=1}^n \delta_{Y_j}\) denote the empirical distribution, and \(\tilde P_n(t) = \frac1n \sum_{j=1}^n e^{-itY_j}\) its Fourier transform.
We assume throughout that \(\tilde g\) is Lebesgue-almost-everywhere non-zero and \(\int f(x)^2\,dx < \infty\).

\begin{definition}
The \emph{smoothness penalized deconvolution} or \emph{SPeD} estimator is defined by the following equivalent characterizations.
\begin{enumerate}[label=(\roman*)]
 \item \(f_n^\alpha = \argmin_{u \in H_{2,m}(\mathbb{R})} \|g*u\|^2 - \frac2n\sum_{j=1}^ng*u(Y_i) + \alpha \|u^{(\nu)}\|^2,\)
 \item \(f_n^\alpha(x) = \frac1n\sum_{i=1}^n \varphi_\alpha(x-Y_i)\), and
 \item \(\tilde f_n^\alpha(t) = \tilde \varphi_\alpha(t) \tilde P_n(t)\),
\end{enumerate}
where \(\tilde \varphi_\alpha(t) = \frac{\overline{\tilde g(t)}}{|\tilde g(t)|^2 + \alpha |t|^{2\nu}}\) and \(\varphi_\alpha(x) = \lim_{r \to \infty} \frac1{2\pi}\int_{-r}^r e^{itx}\tilde \varphi_\alpha(t)\,dt\) is its inverse Fourier transform.
The Sobolev space of functions with square-integrable \(m\)th derivative is denoted by \(H_{2,m}(\mathbb{R})\).
The proof that these are equivalent representations is essentially identical to the proof of Theorem 2(i)-(iii) in \cite{kentSmoothnessPenalizedDeconvolutionSPeD2022}.
\label{def:sped}
\end{definition}

\begin{proposition}\label{prop-combined-fna}
  The following facts about the estimator will be useful later.
  \begin{enumerate}[label=(\roman*)]
  \item \label{prop-combined-fna-iv} \(\int \var(f_n^{\alpha}(x))\,dx = \frac1n \|\varphi_{\alpha}\|^2 - \frac1n \|\varphi_{\alpha}*g*f\|^2\)
  \item \label{prop-combined-fna-l2} \(\|f_n^{\alpha}\|^2 = \frac{1}{n} \|\varphi_{\alpha}\|^2 + \frac{1}{n^2} \sum_{j=1}^n \sum_{k \neq j} \varphi_{\alpha}*\varphi_{\alpha}(x + Y_j - Y_k)\)
  \item \label{prop-combined-fna-El2} \(\mathbb{E}\|f_n^{\alpha}\|^2 = \frac{1}{n} \|\varphi_{\alpha}\|^2 + \frac{n-1}{n} C_{\alpha,f}\). \(C_{\alpha,f}\) does not depend on \(n\).
  \end{enumerate}
  For the following two facts, assume \(\int |f'(x)|^2\,dx < \infty\).
  \begin{enumerate}[label=(\roman*)]  \setcounter{enumi}{3}
  \item\label{prop-combined-fna-normal-rates} (Normal errors) Assume \(g(x) = (2\pi)^{-1} e^{-x^2/2}\). If \(\alpha_n = \kappa \frac{\log n}{n}\), then \(\mise(\alpha_n) = O([\log n]^{-1})\).
  \item\label{prop-combined-fna-cauchy-rates} (Cauchy errors) Assume \(g(x) = \frac{1}{\pi(1+x^2)}\). If \(\alpha_n = \kappa \frac{(\log n)^2}{n}\), then \(\mise(\alpha_n) = O([\log n]^{-2})\).
  \end{enumerate}
\end{proposition}
\begin{proof}[Proof of Proposition~\ref{prop-combined-fna}.]
\textbf{Part~\ref{prop-combined-fna-iv}:} 
  First, since \(f_n^{\alpha} = \frac1n \sum_{j=1}^n\varphi_{\alpha}(x-Y_j)\) is a sum of independent summands, \(\int \var(f_n^{\alpha}(x))\,dx = \frac1n \int \var(\varphi_{\alpha}(x-Y_1))\,dx\).
  Writing the variance in terms of the moments,
  \[\begin{aligned}
   \int \var(\varphi_{\alpha}(x-Y_1))\,dx &= \int \mathbb{E}[\varphi_{\alpha}(x-Y_1)^2]\,dx - \int \mathbb{E}[\varphi_{\alpha}(x-Y_1)]^2\,dx\\
  \end{aligned}\]
For the first,
\[\begin{aligned}
 \int \mathbb{E}[\varphi_{\alpha}(x-Y_1)^2]\,dx &= \int \int \varphi_{\alpha}(x-y)^2h(y)\,dy\,dx\\
                                              &= \int \int \varphi_{\alpha}(x-y)^2\,dx\,h(y)\,dy\\
                                              &= \int \|\varphi_{\alpha}\|^2\,h(y)\,dy\\
                                              &= \|\varphi_{\alpha}\|^2.
\end{aligned}\]
For the second, note that \(\mathbb{E}[\varphi(x-Y_1)] = \varphi_{\alpha}*h(x) = \varphi_{\alpha}*g*f(x)\), so \(\int \mathbb{E}[\varphi_{\alpha}(x-Y_1)]^2\,dx = \|\varphi_{\alpha}*g*f\|^2\).
Combining these gives the desired result.

\noindent\textbf{Part~\ref{prop-combined-fna-l2}:} 
  \[\begin{aligned}
   \|f_n^{\alpha}\|^2 &= \int  \left( \frac{1}{n}\sum_{j=1}^n \varphi_{\alpha}(x-Y_j) \right)^2\,dx\\
                    &= \frac1{n^2} \sum_{j=1}^n \sum_{k=1}^n \int  \varphi_{\alpha}(x-Y_j) \varphi_{\alpha}(x-Y_k)\,dx\\
                    &= \frac1{n^2} \sum_{j=1}^n \int \varphi_{\alpha}(x-Y_j)^2\,dx + \frac1{n^2} \sum_{j=1}^n\sum_{k\neq j} \int  \varphi_{\alpha}(x-Y_j) \varphi_{\alpha}(x-Y_k)\,dx\\
                    &= \frac1{n^2} \sum_{j=1}^n \int \varphi_{\alpha}(x)^2\,dx + \frac1{n^2} \sum_{j=1}^n\sum_{k\neq j} \int  \varphi_{\alpha}(x) \varphi_{\alpha}(x - (Y_k-Y_j))\,dx\\
                    &= \frac1n \|\varphi_{\alpha}\|^2 + \frac1{n^2} \sum_{j=1}^n\sum_{k\neq j} \varphi_{\alpha}*\varphi_{\alpha}(Y_k-Y_j)
  \end{aligned}\]

\noindent\textbf{Part~\ref{prop-combined-fna-El2}:} Here, we just take the expectation of the expression in Part~\ref{prop-combined-fna-l2}.
The only random term is in the summand and its expectation \(C_{\alpha,f}\) is identical for all \(j\) and \(k\), so the value of the double sum is \(n(n-1)C_{\alpha,f}\).

\noindent \textbf{Parts~\ref{prop-combined-fna-normal-rates}-~\ref{prop-combined-fna-cauchy-rates}:} 
Define \(f^{\alpha}(x) = \varphi_{\alpha}*h(x)\) as in \cite{kentSmoothnessPenalizedDeconvolutionSPeD2022}; this is exactly the expected value of the estimator in Definition~\ref{def:sped}.
We have \(\mise(\alpha) = \int (f^{\alpha} - f)^2 + \int \var(f_n^{\alpha})\).
By part~\ref{prop-combined-fna-iv}, \(\int \var(f_n^{\alpha}) \leq \frac{1}{n}\|\varphi_{\alpha}\|^2\), which can be bounded in the following way for \(\alpha < M\).
By dropping \(\alpha|t|^{2\nu}\) in the denominator, it can be seen that \(|\tilde \varphi_\alpha(t)|^2 \leq |\tilde g(t)|^{-2}\).
Applying the inequality \(a + b \geq 2\sqrt{ab}\) to the denominator yields \(|\tilde \varphi_\alpha(t)|^2 \leq \frac12 \alpha^{-\frac12}|t|^{-\nu}\).
Then \(\int |\tilde \varphi_\alpha(t)|^2\,dt \leq \frac1\alpha \int \min \{ M|\tilde g(t)|^{-2}, \frac14|t|^{-2\nu}\}\,dt\).
The integral is finite because \(|\tilde g(t)|^{-2}\) is integrable near zero, while \(|t|^{-2m}\) is integrable away from zero.
Thus \(\int \var(f_n^{\alpha}) \leq \frac{C}{n\alpha}\), with \(C\) depending on \(g\), \(\nu\), and \(M\) but not on~\(f\).

For normal errors, the integrated squared bias is bounded in \cite{kentSmoothnessPenalizedDeconvolutionSPeD2022} Lemma~5(i), so that \(\int (f^{\alpha} - f)^2 \leq C_1[\nu \log(\nu^{-1}\alpha^{-\frac{1}{\nu}})]^{-1}\).
Combining this with the bound for the variance and \(\alpha_n = \kappa \frac{\log n}{n}\) we have
\[
\begin{aligned}
  \mise(\alpha_n) &\leq C_1[\nu \log(\nu^{-1}\alpha^{-\frac{1}{\nu}})]^{-1} + C_2(n\alpha_n)^{-1}\\
                  &= C_1 \nu^{-1} [-\log \nu - \nu^{-1} \left( \log \kappa + \log \log n - \log n \right) ]^{-1} + C_3(\log n)^{-1}\\
                  &\sim C_4 \log n
\end{aligned}
\]

For Cauchy errors, by \cite{kentSmoothnessPenalizedDeconvolutionSPeD2022} Lemma~5(ii), \(\int (f^{\alpha} - f)^2 \leq C_1[\nu \log(\nu^{-1}\alpha^{-\frac{1}{2\nu}})]^{-2}\).
Combining this with the bound for the variance and \(\alpha_n = \kappa \frac{(\log n)^2}{n}\) we have
\[
\begin{aligned}
  \mise(\alpha_n) &\leq C_1\nu^{-2} [\log(\nu^{-1}\alpha_n^{-\frac{1}{2\nu}})]^{-2} + C_2(n\alpha_n)^{-1}\\
                  &= C_1\nu^{-2} [\log \nu + (2\nu)^{-1} \left( \log \kappa + 2 \log \log n - \log n \right)]^{-2} + C_3(\log n)^{-2}\\
                  &\sim C_4(\log n)^{-2}
\end{aligned}
\]

\end{proof}

Note that \(0 \leq \|\varphi_{\alpha}*g*f\|^2 \leq \|f\|^2\) for all \(\alpha\), so Proposition~\ref{prop-combined-fna}~\ref{prop-combined-fna-iv} implies that \(\int \var(f_n^{\alpha}(x))\,dx = \frac1n \|\varphi_{\alpha}\|^2 + O(\frac{1}{n})\).

\section{Cross-Validation of the Error}
\label{sec:criterion}

Cross-validation, proposed in \cite{rudemo_empirical_1982}, is a data-based method for estimating a constant shift of the integrated squared error \(\|\hat f_n - f\|^2\) in an estimate \(\hat f_n\).
The cross-validation criterion has the same expectation as \(\|\hat f_n - f\|^2 - \|f\|^2 = \|\hat f_n\|^2 - 2\int \hat f_n f\).
Since the constant shift \(\|f\|^2\) does not depend on the estimate \(\hat f_n\), this criterion can be used in place of \(\|\hat f_n - f\|^2\) for choosing among possible \(\hat f_n\), typically by setting a tuning parameter.
In the following section we derive the CV criterion for SPeD.

Following that, we introduce the related SCV criterion, which is estimates the error in the estimator fit to a smaller sample size \(m \ll n\); in this context, the particular integrated error \(\|\hat f_m - f\|^2\) makes less sense, so we will focus on the mean integrated squared error (MISE).
Both the CV and SCV criteria can be regarded as an unbiased estimator of a vertical shift of the MISE.
\begin{definition}\label{def-r}
For the SPeD estimator applied to a sample of size \(n\), define
\[\begin{aligned}
  R(\alpha,n) &= \mathbb{E}\|f_n^{\alpha} - f\|^2 - \|f\|^2 \\
              &= \mathbb{E}\|f_n^{\alpha}\|^2 - 2\int \mathbb{E}f_n^{\alpha} f
\end{aligned}\]
\end{definition}

\subsection{Cross-Validation in Density Deconvolution}
\label{sec:criterion:cv}
Least-squares cross-validation for classical (error-free) density estimation can be motivated by noticing that in the integrated squared error of an estimate \(\hat f\), i.e. \(\|\hat f - f\|^2 = \|\hat f\|^2 - 2\int \hat f f + \|f\|^2\), the first term is known and the last term is constant, so we need only estimate the second term.
In classical density estimation, the data \(X_i\) have pdf \(f\), and contemplating the integral \(\int \hat f f\) we might attempt to approximate this by a sample average \(\frac1n \sum \hat f(X_i)\).
This turns out to differ in expectation from \(\int \hat f f\) because \(\hat f\) depends on \(X_i\); this is corrected by using \(\frac1n \sum_{i=1}^n\hat f_{-i}(X_i)\), where \(\hat f_{-i}\) refers to the estimator fit on the data excluding observation \(i\).
The resulting expression \(\|\hat f\|^2 - \frac2n \sum_{i=1}^n\hat f_{-i}(X_i)\) is the (unbiased) cross-validation estimate of the error \citep{scott_biased_1987}.

If we are estimating a density \(f\), but observe data with a different density \(h\) (as we are in the deconvolution problem), this line of thinking grinds to a halt when we think of approximating \(\int \hat f f\) by a sample average -- we have no samples with density \(f\) that we can use to that end.
But we can extend the idea in the following way.
Assume the data \(Y_i\) have density \(h = \mathcal{C}f\) for some linear operator \(\mathcal{C}\) with adjoint \(\mathcal{D}\) (satisfying \(\int u[\mathcal{C}v] = \int [\mathcal{D}u]v\)); assume also that \(\mathcal{D}\) is injective and \(\hat f\) is in the range of \(\mathcal{D}\).
Then we manipulate that integral like so: \(\int \hat f f = \int [\mathcal{D}\mathcal{D}^{-1}\hat f]f = \int [\mathcal{D}^{-1}\hat f][\mathcal{C}f] = \int [\mathcal{D}^{-1}\hat f]h\).
Now we can move forward and approximate \(\int [\mathcal{D}^{-1}\hat f]h\) by the sample average \(\frac1n \sum_{i=1}^n \mathcal{D}^{-1}\hat f_{-i}(Y_i)\).

In deconvolution, we have \(h(y) = \mathcal{C}f(y) = \int f(t)g(y-t)\,dt\), in which case \(\mathcal{D}\) is defined by \(\mathcal{D}u(y) = \int g(t-x)u(t)\,dt\).
The Fourier transform of \(\mathcal{D}u\) is \(\overline{\tilde g}\tilde u\), so our assumption that \(\tilde g\) is a.e. non-zero guarantees that \(\mathcal{D}\) is injective.
For the SPeD estimator, we can see that \(f_n^{\alpha}\) is in the range of \(\mathcal{D}\) by observing that \(\tilde f_n^{\alpha} = \overline{\tilde g} \tilde \varphi_{\alpha} \tilde P_n/\overline{\tilde g}\).
Now \(\tilde u = \tilde \varphi_{\alpha} \tilde P_n/\overline{\tilde g}\) is square-integrable, so its Fourier inverse \(u\) exists; thus \(f_n^{\alpha} = \mathcal{D}u\).
The operator \(\mathcal{D}^{-1}\) can be defined explicitly through Fourier transforms.
\begin{definition}
  \(\mathcal{D}^{-1}u(x) = \frac1{2\pi} \int e^{itx} \tilde u(t)/\overline{\tilde g(t)}\,dt\) for \(u\) in the range of \(\mathcal{D}\).
\end{definition}

Finally, we can state the cross-validation estimate of the error.
It may be useful to note that the same form is obtained if one follows the approach of \cite{youndjeOptimalBandwidthSelection2008}; there, the authors substitute an unbiased estimator of \(|\tilde h(t)|^2\) in \(R(\alpha,n)\) after taking Fourier transforms.
\begin{definition}
  The \emph{cross-validation criterion} is
  \begin{equation}
  \label{eq:cv}
  \cv_n(\alpha) = \|f_n^{\alpha}\|^2 - \frac{2}{n}\sum_{j = 1}^n \mathcal{D}^{-1}f_{n,-j}^{\alpha}(Y_j)
  \end{equation}
  \label{def-cv-criterion}
\end{definition}

\subsection{Stabilized Cross-Validation}
\label{sec:criterion:scv}

The cross-validation criterion performs poorly in practice.
The value of \(\cv_n(\alpha)\) is very unstable near \(\alpha_n = \argmin R(\alpha,n)\), and the minimizer of \(\cv_n(\alpha)\) is often much too small; furthermore, realizations of \(\cv_n(\alpha)\) are often monotonically increasing in an extremely large interval around \(\alpha_0\), and in these cases the minimum is determined in practice by the arbitrary left endpoint.
To remedy this, we introduce a related criterion \(\scv_n(\alpha,m)\) that estimates \(R(\alpha,m)\) with \(m\) possibly different from the sample size \(n\); when \(m \ll n\), \(\scv_n(\alpha,m)\) does not suffer from the deficiencies of \(\cv_n(\alpha)\).
Notice that \(\scv_n(\alpha,n) = \cv_n(\alpha)\).
\begin{definition}
  The \emph{stabilized cross-validation criterion} is
  \begin{equation}
  \label{eq:scv}
  \scv_n(\alpha,m) = \frac{(m-1)n}{m(n-1)} \|f_n^{\alpha}\|^2 + \frac{n-m}{m(n-1)} \|\varphi_{\alpha}\|^2 - \frac{2}{n}\sum_{j = 1}^n \mathcal{D}^{-1}f_{n,-j}^{\alpha}(Y_j) 
  \end{equation}
  \label{def-scv-criterion}
\end{definition}
Just as \(\cv_n(\alpha)\) is an unbiased estimate of \(R(\alpha,n)\), this stabilized criterion \(\scv_n(\alpha,m)\) is unbiased estimate of the error for a sample of size~\(m\).
\begin{proposition}\label{prop-ER}
  \(\mathbb{E}\scv_n(\alpha,m) = R(\alpha,m)\)
\end{proposition}
\begin{proof}[Proof of Proposition~\ref{prop-ER}]
  First, note that \(\mathbb{E}\int f_n^{\alpha} f = \int [\varphi_{\alpha}*h] f\).
  Now, consider the sum: \(\mathbb{E} \frac1n \sum_j \mathcal{D}^{-1}f_{n,-j}^{\alpha}(Y_j) = \mathbb{E}\mathcal{D}^{-1}\varphi_{\alpha}(Y_2-Y_1)\).
  Recalling that \(h=\mathcal{C}f\), we write \(\mathbb{E}\mathcal{D}^{-1}\varphi_{\alpha}(Y_2-Y_1) = \int h(t) \int [\mathcal{D}^{-1}\varphi_{\alpha}](x-t) \mathcal{C}f(x)\,dx\,dt = \int h(t) \int \varphi_{\alpha}(x-t) f(x)\,dx\,dt = \int [\varphi_{\alpha}*h]f\), so that \(\mathbb{E} \frac1n \sum_j \mathcal{D}^{-1}f_{n,-j}^{\alpha}(Y_j) = \mathbb{E} \int f_n^{\alpha} f = \mathbb{E} \int f_m^{\alpha} f\) (since this quantity does not depend on \(n\)), bringing us to \(\mathbb{E}\hat R(\alpha,m) = \frac{(m-1)n}{m(n-1)}\mathbb{E}\|f_n^{\alpha}\|^2 + \frac{n-m}{m(n-1)} \|\varphi_{\alpha}\|^2 - 2\mathbb{E}\int f_m^{\alpha} f \).
  Finally, use Proposition~\ref{prop-combined-fna}~\ref{prop-combined-fna-El2} twice (once with \(n\) and once with \(m\)) to find that \(\frac{(m-1)n}{m(n-1)}\mathbb{E}\|f_n^{\alpha}\|^2 + \frac{n-m}{m(n-1)}\|\varphi_{\alpha}\|^2 = \frac{1}{m}\|\varphi_{\alpha}\|^2 + \frac{m-1}{m}C(\alpha) = \mathbb{E}\|f_m^{\alpha}\|^2\), where \(C(\alpha)\) is the same constant as in Proposition~\ref{prop-combined-fna}~\ref{prop-combined-fna-El2}.
  Bringing this all together, we have \(\mathbb{E}\hat R(\alpha,m) = \mathbb{E}\|f_m^{\alpha}\|^2 - 2\mathbb{E}\int f_m^{\alpha} f = R(\alpha,m)\), as needed.
\end{proof}

\begin{figure}
  \centering
  \includegraphics{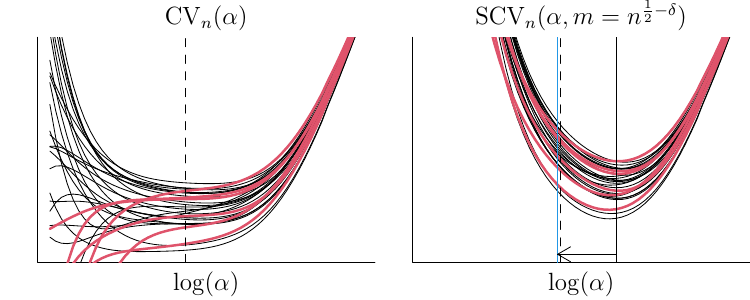}
  \caption{Thirty realizations of the \(\cv\) and \(\scv\) functions under Density 2 (see Section~\ref{sec:simulations}) with Gaussian errors, \(n = 1000\), and 10\% measurement error; \(\scv\) is computed with \(m = n^{\frac{1}{2} - \delta}\) with \(\delta = 10^{-3}\). Realizations plotted in red have \(\cv(\alpha)\) monotonically increasing over the whole domain; the \(\scv_n(\alpha,m)\) functions computed from those same realizations are plotted in red. The dashed and vertical lines are \(\alpha_n^{*} = \argmin R(\alpha,n)\) and \(\alpha_m^{*}\), respectively. The blue vertical line is \((b_n/b_m)\alpha_m^{*}\).}
  \label{fig:realizations}
\end{figure}

The obvious objection here is: ``OK, maybe you've given me a nice way to choose \(\alpha\) when my sample size is \(m\). But my sample size is \(n\)!''
To make the leap from the minimum of \(R(\alpha,m)\) to the minimum of \(R(\alpha,n)\), we take an asymptotic point of view inspired by \cite{hallUsingBootstrapEstimate1990}.
Suppose that the sequence \(\alpha_n^{*} = \argmin_{\alpha \in \mathcal{A}_n} R(\alpha,n)\) satisfies \(\alpha_n^{*} \sim \kappa^{*} b_n\) for some known sequence \(b_n\). 
Then we have that \(\kappa^{*} \sim \alpha_m^{*}/b_m\), so \(\alpha_n^{*} \sim (b_n/b_m) \alpha_m^{*}\).
Now we estimate \(\alpha_m^{*}\) by choosing the \(\alpha_m'\) that minimizes \(\scv_n(\alpha,m)\), and approximate \(\alpha_n^{*}\) by \(\hat \alpha_n = (b_n/b_m) \alpha_m'\).
For example, when the measurement errors are Gaussian and \(f\) has a square-integrable first derivative, \(b_n = n^{-1}\log n\) (Proposition~\ref{prop-combined-fna}~\ref{prop-combined-fna-normal-rates}); in that case we obtain \(\alpha_m'\) by minimizing \(\scv_n(\alpha,m)\), and use \(\hat \alpha_n = \frac{m\log n}{n \log m} \alpha_m'\).
The intuition behind SCV is pictured in Figure~\ref{fig:realizations}.
The procedure can be summarized as follows.
\begin{definition}
  Suppose that \(\alpha_n^{*} = \argmin_{\alpha \in \mathcal{A}_n} R(\alpha,n)\) satisfies \(\alpha_n^{*} \sim \kappa^{*} b_n\) for known sequence \(b_n\), and choose \(m \ll n\).
  The \emph{stabilized cross-validation procedure} (SCV procedure) is defined by
  \begin{enumerate}
  \item Compute \(\alpha_m' = \argmin_{\alpha \in \mathcal{A}_n} \scv_n(\alpha,m)\).
  \item Set \(\hat \alpha_n = (b_n/b_m) \alpha_m'\).
  \item Use \(f_n^{\hat \alpha_n}\) to estimate \(f\).
  \end{enumerate}
  \label{def-scv-procedure}
\end{definition}

The sequence \(b_n\) depends on smoothness assumptions for the target density \(f\); for instance with normal errors, \(b_n = n^{-1}\log n\) attains optimal rates for \(f\) with one square-integrable derivative.
But if we restrict to densities with \(k\) square-integrable derivatives, we should use \(b_n = n^{-1}(\log n)^k\), letting the smoothness parameter decay more slowly \citep{kentSmoothnessPenalizedDeconvolutionSPeD2022}.
In Definition~\ref{def-scv-procedure}, we simply made a minimal smoothness assumption.
We also propose an adaptive procedure for normal errors, in which we choose \(k\) from the data, motived by the fact that \(\alpha_n^{*} \approx \kappa n^{-1} (\log n)^k\), so \(\log \alpha_n^{*} \approx \log \kappa + k\log\log n - \log n\).
\begin{definition}
  Choose \(m_i \ll n\), \(i = 1, \dots, n_m\) and \(m \ll n\).
  The \emph{adaptive stabilized cross-validation procedure} (ASCV procedure) for normal errors is defined by
  \begin{enumerate}
  \item Compute \(\alpha_{m_i}' = \argmin_{\alpha \in \mathcal{A}_n} \scv_n(\alpha,m_i)\) for \(i = 1, \dots n_m\).
  \item Fit the regression model \(\log \alpha_{m_i}' = \beta_0 + \beta_1\log \log m_i - \log m_i + \varepsilon_i\), and set \(\hat b_n = n^{-1} (\log n)^{\hat \beta_1}\).
  \item Set \(\hat \alpha_n = (\hat b_n/\hat b_m) \alpha_m'\).
  \item Use \(f_n^{\hat \alpha_n}\) to estimate \(f\).
  \end{enumerate}
  \label{def-ascv-procedure}
\end{definition}

\begin{remark}[Choosing \(m\)]\label{remark:choosing-m}
  Definitions~\ref{def-scv-procedure}~and~\ref{def-ascv-procedure} both require choice of \(m \ll n\).
  The SCV procedures both involve a statistical estimate and an asymptotic approximation: the value \(\alpha_m'\) is a statistical estimate of the optimal \(\alpha_m^{*}\), and we then use the asymptotic approximation \(\alpha_n^{*} \approx (b_n/b_m)\alpha_m^{*}\).
  The smaller \(m\) is with respect to \(n\), the more well-behaved are the estimates \(\alpha_m'\); on the other hand, the asymptotic approximation \(\alpha_n^{*} \approx (b_n/b_m) \alpha_m^{*}\) relies on both \(n\) and \(m\) being large.
  So the trade-off is that when \(m\) is very small, the estimate \(\alpha_m' \approx \alpha_m^{*}\) is quite good, but the asymptotic approximation \(\alpha_n^{*} \approx (b_n/b_m)\alpha_m^{*}\) may be poor, while when \(m\) is closer to \(n\), the situation is reversed, and the \(\alpha_m'\) is poorly behaved (cf. Figure~\ref{fig:realizations}) while the approximation is quite good.

  In Theorem~\ref{theorem-contsupdiff_as}, we prove uniform consistency for the \(\scv_n(\alpha,m)\) criterion estimating \(R(\alpha,m)\) under Assumption~\ref{assumption:power} relating \(m\) to the rate of convergence \(r_n = R(\alpha_n^{*},n)\) of the optimal \(\alpha_n^{*} \sim \kappa b_n\).
  For normal errors, we show in Example~\ref{example:normal} that the assumption is equivalent to \(m = n^{\frac{1-\delta}{2}}\) with \(0 < \delta < 1\).
  We recommend choosing \(m\) large, but still satisfying this condition; in all the simulations here, we use \(m = n^{\frac{1}{2} - 10^{-3}}\).
  In practice, \(m = \sqrt{n}\) seems to work just as well, but we did notice degradation of performance for larger powers like \(m = n^{\frac{3}{5}}\) and above.

  When choosing \(m_i\) in the ASCV procedure, we recommend using equally-spaced values close to, but not exceeding \(m\); we found that \(m_1 = n^{\frac{1}{2} - \frac{1}{20}}\) through \(m_{n_m} = n^{\frac{1}{2} - 10^{-3}}\) worked well. 
\end{remark}

\subsection{Relationship between the CV and SCV Criteria}

In view of Proposition~\ref{prop-combined-fna}~\ref{prop-combined-fna-iv}, define
\[
  \ive(\alpha,n) = \frac1n \|\varphi_{\alpha}\|^2.
\]
This approximation to the integrated variance satisfies \(\ive(\alpha,n) = \int \var(f_n^{\alpha}(x))\,dx + O(\frac1n)\).
For reasons we will explain momentarily, define
\[
  \ibsh(\alpha) = \frac1{n^2} \sum_{j=1}^n \sum_{k \neq j} \varphi_{\alpha}*\varphi_{\alpha}(Y_j-Y_k) - \frac2n \sum_{j=1}^n \mathcal{D}^{-1}f_{n,-j}^{\alpha}(Y_j).
\]
By Proposition~\ref{prop-combined-fna}~\ref{prop-combined-fna-l2}, we may write
\[
  \cv_n(\alpha) = \ive(\alpha,n) + \ibsh(\alpha).
\]
By Proposition~\ref{prop-ER} with \(m=n\),
\[\begin{aligned}
  \mathbb{E}\cv_n(\alpha) &= \|f_n^{\alpha} - f\|^2 - \|f\|^2\\
                        &= \int \var(f_n^{\alpha}(x))\,dx + \int (\mathbb{E}f_n^{\alpha}(x)-f(x))^2\,dx - \|f\|^2,
\end{aligned}\]
so
\[
  \mathbb{E}\ibsh(\alpha) = \int (\mathbb{E}f_n^{\alpha}(x)-f(x))^2\,dx - \|f\|^2 + O(\frac1n).
\]
Thus \(\ibsh(\alpha)\) estimates a vertical shift of the integrated squared bias and furthermore is the only random term in \(\cv_n(\alpha)\).
Exactly this observation is made for cross-validation in classical density estimation in \cite{hall_smoothed_1992}.

Using Proposition~\ref{prop-combined-fna}~\ref{prop-combined-fna-l2} again, we can write
\[\begin{aligned}
  \scv_n(\alpha,m) &= \ive(\alpha,m) + \frac{m-1}{mn(n-1)} \sum_{j=1}^n \sum_{k \neq j} \varphi_{\alpha}*\varphi_{\alpha}(Y_j-Y_k)\\
                   &\qquad\qquad\qquad\qquad- \frac2n \sum_{j=1}^n \mathcal{D}^{-1}f_{n,-j}^{\alpha}(Y_j)\\
                 &\approx \frac{n}{m}\ive(\alpha,n) + \ibsh(\alpha),
\end{aligned}\]
showing that to a good approximation, the \(\scv_n(\alpha,m)\) and \(\cv_n(\alpha)\) criteria contain the same estimates of the bias and variance, with the variance rescaled for the smaller sample size.
This shows why the \(\scv_n(\alpha,m)\) is better behaved: it is because \(\var(\ibsh(\alpha))\) decreases as \(\alpha\) increases, and when \(m \ll n\), the minimum of \(R(\alpha,m)\) (and \(\scv_n(\alpha,m)\)) occurs at larger \(\alpha\).

\section{Theory for \(\scv_n(\alpha,m)\)}
\label{sec:theory}

The criterion \(\scv_n(\alpha,m)\) can be expressed as a constant (depending on \(\alpha\) and \(m\)) plus a U-statistic.
From this, it can be shown that \(\scv_n(\alpha,m)\) is the minimum variance unbiased estimator of \(R(\alpha,m)\).
\begin{proposition}\label{prop-ustat}
  \begin{equation}
      \scv_n(\alpha,m) = \frac{1}{m} \|\varphi_{\alpha}\|^2 + \frac{1}{\binom{n}{2}}\sum_{1 \leq j < k \leq n} \theta_\alpha(Y_j - Y_k),
  \end{equation}
  where \(\theta_\alpha(x) = \frac{m-1}{m}\varphi_{\alpha}*\varphi_{\alpha}(x) - 2\mathcal{D}^{-1}\varphi_{\alpha}(x)\).
  Furthermore, for \(k < 2\nu - 1\) the kernel \(\theta_\alpha(x)\) has the bounds
  \begin{equation}
    \sup_x |\theta_\alpha^{(k)}(x)| \leq \frac{C_{M,k}}{\alpha}
  \end{equation}
  for all \(\alpha < M\), where the constant \(C_{M,k}\) depends on \(M\) and \(k\).
\end{proposition}

\begin{proposition}\label{prop-minvar}
  If \(\hat R_n(\alpha,m) = \hat R(\alpha,m;Y_1, \dots, Y_n)\) is an unbiased estimator of \(R(\alpha,m)\), then \(\var(\hat R_n(\alpha,m)) \geq \var(\scv_n(\alpha,m))\)
\end{proposition}

Now we proceed to some consistency results for \(\scv_n(\alpha,m)\) in estimating \(R(\alpha,m)\).
\begin{assumptions*}\label{assumptions}
  Fix \(g\), \(f\), and \(0 < \lambda < \mu\).
  \begin{enumerate}[label=(A\arabic*)]
  \item\label{assumption:rate} There are sequences \(b_n\) and \(r_n\) such that \(\inf_{\lambda b_n \leq \alpha \leq \mu b_n} R(\alpha,n) = r_n\).
  \item\label{assumption:power} The sequence \(m_n\) is chosen such that \((r_{m_n}b_{m_n})^2n = n^{\delta}\) with \(\delta > 0\).
  \end{enumerate}
\end{assumptions*}

Under Assumptions~\ref{assumption:rate} and \ref{assumption:power}, we can show that with probability 1, the stabilized cross-validation criterion converges uniformly to the shifted mean integrated squared error that it is estimating.
\begin{theorem}\label{theorem-contsupdiff_as}
  Assume \ref{assumption:rate} and \ref{assumption:power}. Then with probability 1,
  \[
    \sup_{\lambda b_m \leq \alpha \leq \mu b_m} |\scv_n(\alpha,m) - R(\alpha,m)| = o(r_m).
  \]
\end{theorem}
A simple corollary says that the minimizer of \(\scv_n(\alpha,m)\) asymptotically minimizes the true risk function \(R(\alpha,m)\).
Note that the function \(R(\cdot,m)\) here still refers to the risk of the estimator using a fixed penalty parameter.
\begin{corollary}\label{cor-minimasupdiff}
  Assume \ref{assumption:rate} and \ref{assumption:power}. Then with probability 1,
  \[
    R(\alpha_m',m) = R(\alpha_m^{*},m) + o(r_m),
  \]
  so
  \[
    \frac{R(\alpha_m',m)}{R(\alpha_m^{*},m)} \sim 1.
  \]
\end{corollary}

As a bonus, an observation of \cite{chacon_estimation_2024} leads directly from Corollary~\ref{cor-minimasupdiff} to a strongly consistent estimator of \(\|f\|^2\).
\begin{corollary}\label{cor-f2normest}
  Assume \ref{assumption:rate} and \ref{assumption:power}.
  Then with probability 1,
  \[
    -\scv_n(\alpha_m',m) = \|f\|^2 + O(r_m).
  \]
Since \(r_m \to 0\), it follows that \(-\scv_n(\alpha_m',m) \to \|f\|^2\) with probability 1.
\end{corollary}

\begin{example}[Normal errors]\label{example:normal}
  If \(g\) is a normal density and we assume \(f\) has one square-integrable derivative, then by Proposition~\ref{prop-combined-fna}~\ref{prop-combined-fna-normal-rates} \(\alpha_n = \kappa b_n = \kappa \frac{\log n}{n}\) is optimal, and yields rate of convergence \(r_n = (\log n)^{-1}\). Then \((b_{m} r_{m})^2 = \kappa^2m^{-2}\), so Theorem~\ref{theorem-contsupdiff_as} holds if we use \(m = n^{\frac{1-\delta}2}\) for \(0 < \delta < 1\).
\end{example}

\begin{example}[Cauchy errors]\label{example:cauchy}
 If \(g\) is a Cauchy density and we assume \(f\) has one square-integrable derivative, then by Proposition~\ref{prop-combined-fna}~\ref{prop-combined-fna-cauchy-rates} \(\alpha_n = \kappa b_n = \kappa \frac{(\log n)^2}{n}\) is optimal, and yields rate of convergence \(r_n = (\log n)^{-2}\). Then \((b_mr_m)^2 = \kappa^2m^{-2}\), and again if we use \(m = n^{\frac{1-\delta}2}\) for \(0 < \delta < 1\) then Theorem~\ref{theorem-contsupdiff_as} holds.
\end{example}

\section{Simulations}
\label{sec:simulations}

\newcommand{\ASCVbettertvalue}{28.78}
\newcommand{\ASCVworst}{1}
\newcommand{\ASCVbeatsSURE}{22}
\newcommand{\ASCVbeatsDKDE}{16}
\newcommand{\maxpropASCV}{81}
\newcommand{\maxpropSURE}{670}
\newcommand{\maxpropDKDE}{254}

\begin{figure}
  \centering
  \includegraphics{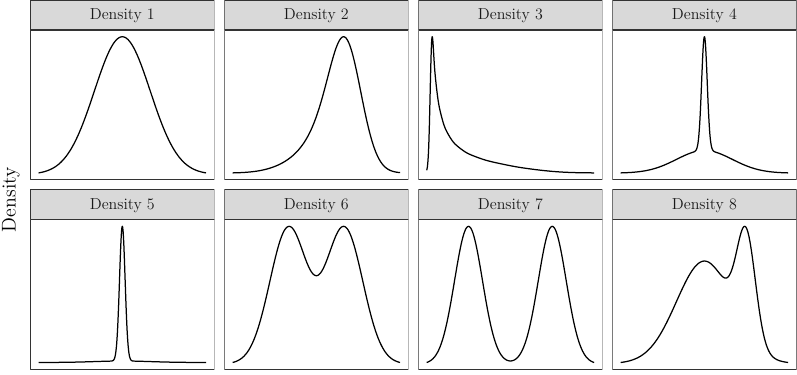}
  \caption{The first eight target densities of \cite{marronExactMeanIntegrated1992} which we test against in simulations.}
  \label{fig:examples}
\end{figure}

To assess the practical performance of the SCV approach, we test it on a battery of densities with a variety of features.
In \cite{marronExactMeanIntegrated1992}, the authors propose fifteen densities which exhibit a variety of challenging features for curve estimation, all represented as mixtures of normal densities.
We use here the first eight, which are the unimodal and bimodal examples, all pictured in Figure~\ref{fig:examples}.
Density 1 is the standard normal; densities 2 \& 3 are slightly and strongly skewed, respectively; density 4 has \(\frac{1}{3}\) of its mass in a sharp peak at the mean; density 5 has \(\frac{9}{10}\) of its mass in a sharp peak at its mean (or from another point of view, \(\frac{1}{10}\) of its mass spread far from the mean); densities 6 and 7 are bimodial with slightly and very separated modes, respectively; density 8 is bimodal with modes of different width.

We aim to estimate each of these target example densities under measurement error.
We simulate \(Y_i = X_i + E_i\), \(i = 1, \dots, n\) with the \(X_i\) drawn independently from the target density \(f\) and \(E_i \sim \N(0,\sigma_E^2)\) with \(\sigma_E^2 = \frac{1}{9}\var(X_1)\), which results in \(\frac{\var(X_i)}{\var(Y_{i})} = \frac{1}{10}\), so that \(10\%\) of the variance in our data is from the measurement error.
Simulations are performed with \(n = 100\), \(500\), and \(1000\).
For each \(f\) and \(n\), we compute \(\alpha_n^{*} = \argmin_{\alpha} \mathbb{E}\|f_n^{\alpha} - f\|^2\), the MISE-optimal smoothing parameter; we will call this the ``optimal'' choice since among fixed choices of \(\alpha\) it has the smallest MISE.
For each iteration, we choose \(\hat \alpha_{\scv}\) by the adaptive SCV (ASCV) procedure described in Definition~\ref{def-ascv-procedure} with \(m = n^{\frac{1}{2} - 10^{-3}}\) and five equally spaced \(m_i\) between \(m_1 = n^{\frac{1}{2} - \frac{1}{20}}\) and  \(m_5 = m\) (see Remark~\ref{remark:choosing-m} for suggestions about choosing \(m\)); we also choose \(\hat \alpha_{\sure}\) by the SURE procedure of \cite{yang_density_2020}.
We fit the estimator \(f_n^{\alpha}\) and compute the ISE \(\|f_n^{\alpha} - f\|^2\) for each of \(\alpha = \alpha_{n}^{*}\), \(\hat \alpha_{\ascv}\), and \(\hat \alpha_{\sure}\); we also fit the classical deconvoluting kernel estimator (DKDE) using the data-based bandwidth choice recommended by \cite{delaigle_practical_2004} and compute its ISE.
For each density \(f\) and sample size \(n\) we estimate the MISE from \(n_{\text{sim}} = 10^4\) iterations.

To compare the ASCV, SURE, and DKDE procedures, we compare the proportional increase in MISE over the optimal \(\alpha_n^{*}\), for example we look at \(100 \times (\mise_{\ascv}/\mise_{\text{optimal}} - 1)\) for the ASCV procedure.
These values are displayed in Table~\ref{table-pmise}.
In \ASCVbeatsSURE{} of the 24 simulation settings, our proposed ASCV procedure outperforms the SURE procedure.
The SPeD estimator with our ASCV procedure outperformed the DKDE with plug-in bandwidth in \ASCVbeatsDKDE{} of the 24 settings.
Another interesting observation is that the ASCV procedure's worst case yields a MISE \maxpropASCV{}\% larger than the optimal choice of parameter; SURE and DKDE have worst cases with MISE \maxpropSURE{}\% and \maxpropDKDE{}\% larger than optimal, respectively.

In assessing the SURE procedure, the authors saw that it occasionally ``selects too small a [penalty parameter],'' and in their Figure~3 they show an example of a catastrophically small penalty parameter, where the estimate has many narrow modes, though the target density is unimodal \citep{yang_density_2020}.
These simulations bear this out--in Figure~\ref{fig:histograms}, the small penalty parameters are apparent in the long left tails of the SURE histograms, and in Figure~\ref{fig:alpha_vs_ise}, we can see that when the SURE estimate has extremely large error, this is typically when it has selected a very small penalty parameter.
The ASCV procedure does not have this same catastrophic failure mode.
It is therefore worth quantifying: how extreme is each procedure's ``poor'' estimates?
We'll consider e.g. \(\ise_{\ascv}/\ise_{\optimal}\), so that an ASCV estimate is poor if its error is much larger than the error of the estimate using \(\alpha_n^{*}\).
We then consider the 99\textsuperscript{th} percentile of this quantity, plotted for three target densities in Figure~\ref{fig:muchworse}; density 2 was the worst case for both SURE and ASCV, density 1 was the worst case for DKDE, and density 5 was the best case for both ASCV and SURE.
The results corroborate the comments of \cite{yang_density_2020}: the SURE procedure's poor estimates may be \emph{much} worse than optimal, even 200 times larger.
On the other hand, the ASCV and DKDE procedure's poor estimates are more pedestrian, around 8-16 times larger.

In Table~\ref{table-pmise}, the negative value for ASCV in target density 5, \(n = 1000\) may arouse suspicion; the fact that this number is negative means that the ASCV choice of tuning parameter resulted in a \emph{lower} MISE than the what we have called the optimal choice of tuning parameter.
First, recall that our notion of optimality involved a fixed, non-random tuning parameter; \(\alpha_n^{*}\) is, for a given \(f\) and \(n\), a fixed value that minimizes the MISE.
There are certainly data-dependent choices of \(\alpha\) which could give lower MISE; for example \(\hat \alpha_n^{*} = \argmin_{\alpha} \|f_n^{\alpha} - f\|^2\) minimizes the error for each realization, and has MISE smaller than \(\alpha_n^{*}\) and all other choices of \(\alpha\).
Since ASCV is data-dependent, it is plausible that it can result in MISE lower than \(\alpha_n^{*}\).
Testing the hypothesis that \(\mise_{\text{ASCV}} < \mise_{\optimal}\) in our simulations results in a \(t\)-statistic of \(t = \ASCVbettertvalue{}\), so we can conclude that in this one simulation setting, the ASCV choice is indeed better than the fixed value \(\alpha_n^{*}\) that minimizes the MISE.

\begin{figure}
  \centering
  \includegraphics{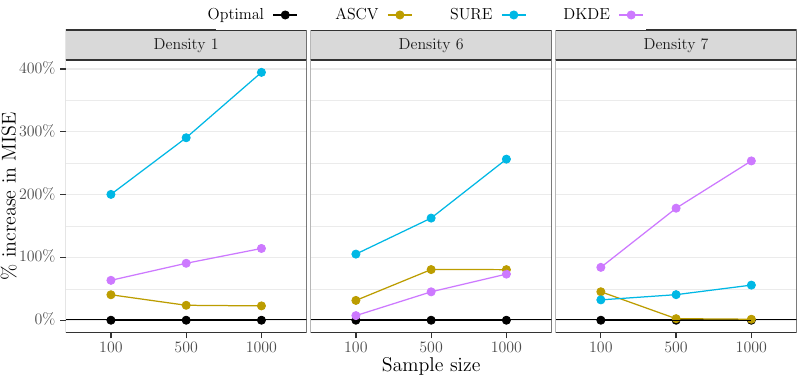}
  \caption{Based on \(n_{\text{sim}} = 10^4\) simulations, the percent increase in MISE over the MISE-optimal SPeD. For example, the ASCV values are \((\text{MISE}_{\text{ASCV}}/\text{MISE}_{\text{Optimal}} - 1)\times 100\%\).}
  \label{fig:MISE}
\end{figure}

\begin{figure}
  \centering
  \includegraphics{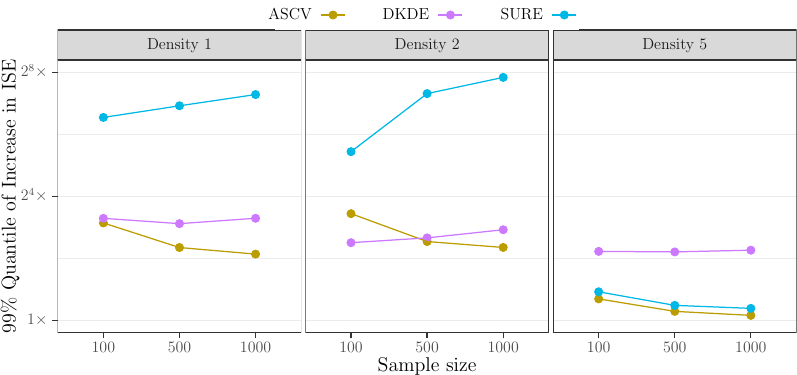}
  \caption{Based on \(n_{\text{sim}} = 10^4\) simulations. The \(0.99\)-quantile of increase in ISE over the MISE-optimal penalty. For example, the ASCV values are the \(0.99\)-quantile of \(\text{ISE}_{\text{ASCV}}/\text{ISE}_{\text{Optimal}}\). In other words: benchmarked against the MISE-optimal parameter, if we have done much worse than usual, how much worse have we done? Note the logarithmic \(y\)-axis.}
  \label{fig:muchworse}
\end{figure}

{\renewcommand{\baselinestretch}{1.05} %
\begin{table}
\centering
\begin{tabular}{llrrrrrrrr}
  & & \multicolumn{8}{c}{Target density}\\\cline{3-10}
  & Method & 1 & 2 & 3 & 4 & 5 & 6 & 7 & 8 \\ 
  \hline
 & ASCV & \textbf{41} & 49 & \textbf{13} & 13 & \textbf{15} & 32 & 45 & 17 \\ 
  \(n = 100\) & SURE & 200 & 145 & 23 & \textbf{8}* & 24 & 105 & \textbf{32}* & 80 \\ 
   & DKDE & 64 & \textbf{43} & 16 & 12 & 106 & \textbf{7} & 84 & \textbf{3} \\ 
   \hline
 & ASCV & \textbf{24} & \textbf{23} & \textbf{9} & \textbf{5} & \textbf{1} & 81 & \textbf{2} & 48 \\ 
  \(n = 500\) & SURE & 291 & 396 & 11 & 9 & 9 & 163 & 41 & 127 \\ 
   & DKDE & 91 & 74 & 27 & 23 & 164 & \textbf{45} & 178 & \textbf{23} \\ 
   \hline
 & ASCV & \textbf{23} & \textbf{23} & \textbf{8} & \textbf{1} & \textbf{-2} & 81 & \textbf{2} & 57 \\ 
  \(n = 1000\) & SURE & 395 & 670 & 11 & 3 & 6 & 256 & 56 & 181 \\ 
   & DKDE & 114 & 100 & 34 & 20 & 186 & \textbf{73} & 254 & \textbf{39} \\ 
   \hline
\end{tabular}
\caption{For the target densities shown in Figure~\ref{fig:examples}, the percent increase in MISE over the MISE-optimal SPeD. For each target density/sample size, the best method (incurring the smallest increase in MISE) is marked in bold. In the two settings where the SURE procedure outperforms the ASCV procedure, the value is marked with an asterisk.}
\label{table-pmise}\end{table}
}

\begin{figure}
  \centering
\includegraphics{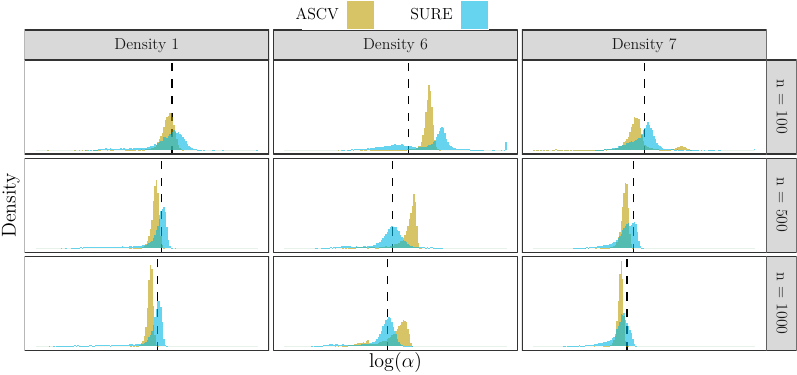}
  \caption{Histograms of \(\alpha_{\ascv}\) and \(\alpha_{\sure}\) in the \(n_{\text{sim}} = 10^4\) simulations. Vertical dashed line marks the MISE-optimal value \(\alpha_n^{*}\) in each setting. The only panel here in which \(\sure\) has lower MISE than \(\ascv\) is Density 7, \(n = 100\).}
  \label{fig:histograms}
\end{figure}

\begin{figure}
  \centering
  \includegraphics{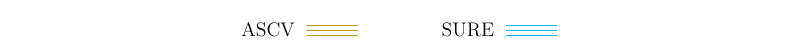}
  \includegraphics{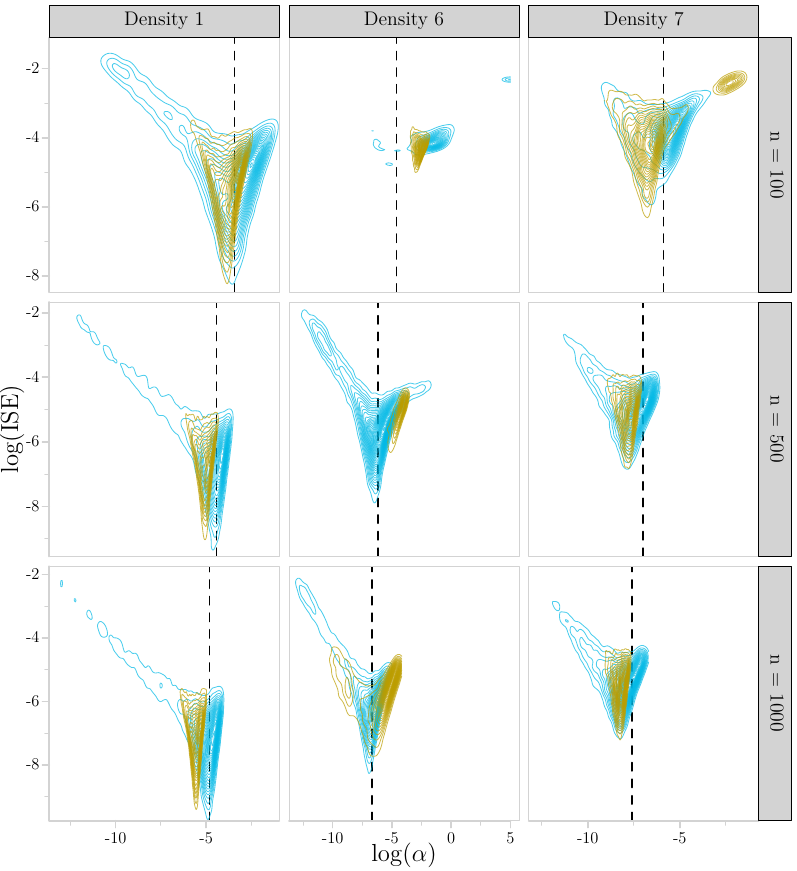}
  \caption{Contour plots of the probability density of \((\log \alpha, \log \ise)\) for ASCV and SURE in the \(n_{\text{sim}} = 10^4\) simulations. Vertical dashed line marks the MISE-optimal value \(\alpha_n^{*}\) in each setting.}
  \label{fig:alpha_vs_ise}
\end{figure}

\section*{Acknowledgements}

The work on this paper would not have been possible without the guidance, support, and encouragement of my PhD advisor David Ruppert.
\par

\section{Appendix}
\label{sec:appendix}

\subsection{Proofs of Results}\label{sec:appendix:proofs}
\begin{proof}[Proof of Proposition~\ref{prop-ustat}]
  Using Proposition~\ref{prop-combined-fna}~\ref{prop-combined-fna-l2} and rearranging,
  \[\begin{aligned}
    \scv_n(\alpha,m) &= \frac{1}{m} \|\varphi_{\alpha}\|^2\\
                     &\qquad+ \frac{2}{\binom{n}{2}} \sum_{1 \leq j < k \leq n} \left[ \frac{m-1}{m}\varphi_{\alpha}*\varphi_{\alpha} - 2\mathcal{D}^{-1}\varphi_{\alpha}\right](Y_j-Y_k) .
  \end{aligned}\]
Now, the Fourier transform of \(\theta_{\alpha}(x)\) is
\[\begin{aligned}
  \tilde \theta_{\alpha}(t) &= \frac{m-1}{m}\frac{|\tilde g(t)|^2}{(|\tilde g(t)|^2 + \alpha |t|^{2\nu})^2} - \frac{2}{|\tilde g(t)|^2 + \alpha |t|^{2\nu}}\\
                           &= \frac{m-1}{m}\frac{|\tilde g(t)|^2}{(|\tilde g(t)|^2 + \alpha |t|^{2\nu})^2} - \frac{2|\tilde g(t)|^2 + 2\alpha |t|^{2\nu}}{(|\tilde g(t)|^2 + \alpha |t|^{2\nu})^2}\\
                           &= -\frac{\frac{m+1}{m}|\tilde g(t)|^2 + 2\alpha |t|^{2\nu}}{(|\tilde g(t)|^2 + \alpha |t|^{2\nu})^2},
\end{aligned}\]
and then \(|\tilde \theta_{\alpha}(t)| \leq \frac{2}{|\tilde g(t)|^2 + \alpha |t|^{2\nu}}\), so that
\[\begin{aligned}
  \theta_{\alpha}^{(k)}(x) &\leq \frac1{2\pi} \int \frac{|t|^k}{|\tilde g(t)|^2 + \alpha |t|^{2\nu}}\,dt\\
                        &= \frac1{\pi} \int_0^{\infty} \frac{t^k}{|\tilde g(t)|^2 + \alpha t^{2\nu}}\,dt.
\end{aligned}\]
Now we integrate this bound near and away from zero.
Choose \(\varepsilon>0\) s.t. \(\tilde g(t) \geq c > 0 \) for all \(0 \leq t \leq \varepsilon\), and bound as
\[\begin{aligned} %
  \theta_{\alpha}^{(k)}(x) &\leq \frac1{\pi} \int_0^{\varepsilon} \frac{t^k}{|\tilde g(t)|^2 + \alpha t^{2\nu}}\,dt + \frac1{\pi} \int_\varepsilon^{\infty} \frac{t^k}{|\tilde g(t)|^2 + \alpha t^{2\nu}}\,dt\\
                         &\leq \frac1{\pi} \int_0^{\varepsilon} \frac{t^k}{|\tilde g(t)|^2 + \alpha t^{2\nu}}\,dt + \frac1{\alpha\pi} \int_\varepsilon^{\infty} t^{k-2\nu}\,dt\\
                         &\leq C_k + C_k'\alpha^{-1}\\
                         &= C_{k,M}\alpha^{-1}\\
\end{aligned}\]
\end{proof}

\begin{proof}[Proof of Proposition~\ref{prop-minvar}]
  The basic idea that follows is similar to Theorem 3 of \cite{lee_u-statistics:_1990}.
  The estimand can be written \(R(\alpha,m) =  \int \int \theta_{\alpha}(y_1-y_2)h(y_1)h(y_2)\,dy_1\,dy_2\), with kernel \(\theta_x(x) = \frac{1}{m} \|\varphi_{\alpha}\|^2 + \frac{m-1}{m} \varphi_{\alpha}*\varphi_{\alpha}(x) - 2\mathcal{D}^{-1}\varphi_{\alpha}(x)\), so \(R(\alpha,m)\) is a functional of degree 2 defined for densities in \(\mathcal{H} = \{h = g*f|f \in \mathcal{F}\}\), with \(\mathcal{F} = \{f \in L_2|\int f = 1, f \geq 0\}\); the U-statistic in Proposition~\ref{prop-ustat} is a symmetric, unbiased estimate formed with this same kernel \(\theta_{\alpha}(x)\).

  To prove the minimum variance assertion, suppose we knew that the symmetric unbiased estimator was unique; i.e. if \(\hat U = U(Y_1, \dots, Y_n)\) is symmetric and \(\mathbb{E}\hat U = R(\alpha,m)\), then \(U(y_1, \dots, y_n) = \frac{1}{\binom{n}{2}}\sum_{1 \leq j < k \leq n} \theta_{\alpha}(y_j - y_k)\) (we will prove this momentarily).
  Let \(\hat R = \hat R(y_1, \dots, y_n)\) be any unbiased estimator of \(R(\alpha,m)\); form a symmetric estimator by averaging \(\hat R(y_1, \dots, y_n)\) over all permutations of the arguments: \(\hat R^{[n]}(y_1, \dots, y_n) = \sum \frac{1}{n!} \hat R(y_{i_1}, y_{i_2}, \dots, y_{i_n})\).
  This is also unbiased.
  Since the symmetric unbiased estimator is unique, \(\hat R^{[n]}(y_1, \dots, y_n) = \scv_n(\alpha,m)\).
  Now apply Cauchy-Schwarz:
  \[\begin{aligned}
   \left( \scv_n(\alpha,m) \right)^2 &= \left( \hat R^{[n]}(y_1, \dots, y_n) \right)^2\\
                                     &= \left( \sum \frac{1}{n!} \hat R(y_{i_1}, y_{i_2}, \dots, y_{i_n}) \right)^2\\
                                     &\leq \sum \left (\frac{1}{n!}\right)^2 \sum \left( \hat R(y_{i_1}, y_{i_2}, \dots, y_{i_n}) \right)^2\\
                                     &= \frac{1}{n!}\sum \left( \hat R(y_{i_1}, y_{i_2}, \dots, y_{i_n}) \right)^2\\
  \end{aligned}\]
Taking expectations,
  \[\begin{aligned}
   \mathbb{E}[\left( \scv_n(\alpha,m) \right)^2] &\leq \frac{1}{n!}\sum \mathbb{E}\left [ \left( \hat R(y_{i_1}, y_{i_2}, \dots, y_{i_n}) \right)^2 \right]\\
                                                 &= \mathbb{E}\left [ \left( \hat R(y_{i_1}, y_{i_2}, \dots, y_{i_n}) \right)^2 \right]
  \end{aligned}\]
Since \(\mathbb{E}\scv_n(\alpha,m) = \mathbb{E} \hat R\), \(\var(\scv_n(\alpha,m)) \leq \var(\hat R)\).

Now we show that the symmetric unbiased estimator is unique.
It suffices to show that for symmetric kernels \(\psi(y_1,y_2, \dots, y_n)\), if
\[
  \int \cdots \int \psi(y_1,\dots,y_n) h(y_1)\dots h(y_n)\,dy_1\dots dy_n = 0
\]
for all \(h \in \mathcal{H}\), then \(\psi(y_1,\dots,y_n) = 0\) for almost all \(y_1,\dots,y_n \in \mathbb{R}\).
Recall that if \(h \in \mathcal{H}\), then \(h = g*f\) for \(f \in \mathcal{F}\).
So we may define \(K(x_1, \dots, x_n) = \int \dots \int \psi(x_1+t_1, \dots, x_n + t_n)g(t_1)\dots g(t_n)\,dt_1\dots dt_n\) and write instead
\[
  \int \cdots \int K(x_1,\dots,x_n) f(x_1)\dots f(x_n)\,dx_1\dots dx_n = 0
\]
for all \(f \in \mathcal{F}\).
Note that \(K\) is symmetric due to the symmetry of \(\psi\).
We will be finished if we show that \(K(x_1, \dots, x_n) = 0\) for almost all \(x_1, \dots, x_n\) for the following reason.
A function vanishes a.e. if and only if its Fourier transform vanishes a.e.; so if we show that \(K(x_1, \dots, x_n)\) vanishes a.e., then so too does its Fourier transform \(\tilde K\).
Now, the Fourier transform of \(K\) is \(\tilde K(s_1, \dots, s_n) = \tilde \psi(s_1, \dots, s_n) \overline{\tilde g(s_1) \dots \tilde g(s_n)}\).
Since this is zero a.e., and \(\tilde g(\cdot)\) is non-zero a.e., we conclude that \(\tilde \psi(s_1, \dots, s_{n}) = 0\) a.e., which ultimately implies that \(\psi(x_1, \dots, x_n) = 0\) for almost all \(x_1, \dots x_n\), which is what we seek.

Consider the symmetric multilinear form
\[
  M(u_1, \dots, u_n) = \int \cdots \int K(x_1,\dots,x_n) u_1(x_1)\dots u_n(x_n)\,dx_1\dots dx_n.
\]
By Lemma~\ref{lemma-multilinear} and the hypothesis that \(M(f,f,\dots,f) = 0\) for all \(f \in \mathcal{F}\), we know that \(M(u_1, u_2, \dots, u_n) = 0\) whenever \(u_1, \dots, u_n \in \vecspan(\mathcal{F})\).
For any Borel set \(B\), the indicator \(\bone_B \in \vecspan(\mathcal{F})\), so the indicator of any product of \(n\) Borel sets can be written as \(u_1(x_1)\dots u_n(x_n)\) with \(u_1, \dots, u_n \in \vecspan(\mathcal{F})\).
Thus \(\int_B K(x_1, \dots, x_n)\,dx_1\dots dx_n = 0\) for any product \(B\) of Borel sets, so \(K(x_1, \dots, x_n) = 0\) almost everywhere.
\end{proof}

\begin{proof}[Proof of Theorem~\ref{theorem-contsupdiff_as}.]
  Let \(\mathcal{A}_n = \{\alpha_i: \alpha_i = \lambda b_m + \frac{i-1}{k_n-1}(\mu - \lambda)b_m\}\), with \(k_n = \lceil n^p \rceil\), \(p > 1\).
  Then Lemma~\ref{lemma-finitesupdiff} holds for \(\alpha_i \in \mathcal{A}_n\), and in fact the probabilities are summable, so it holds almost surely, so with probability 1,
  \[
    \sup_{\alpha \in \mathcal{A}_n} |\scv_n(\alpha,m) - R(\alpha,m)| = o(r_m).
  \]
  Fix an arbitrary \(\alpha' \in [\lambda b_m, \mu b_m]\).
  We want to show that
  \[
    |\scv_n(\alpha',m) - R(\alpha',m)| \leq C_n r_m
  \]
  for some sequence \(C_n \to 0\) which does not depend on \(\alpha'\).
  Let \(\alpha_j \in \mathcal{A}_n\) be the element nearest to \(\alpha'\), and upper bound the difference as
  \begin{equation}
  \label{eq-addsub}
  \begin{aligned}
    |\scv_n(&\alpha',m) - R(\alpha',m)| \\
            &= |\scv_n(\alpha',m) - \scv_n(\alpha_j,m) \\
            & \qquad\qquad - R(\alpha',m) + R(\alpha_j,m) + \scv_n(\alpha_j,m) - R(\alpha_j,m)|\\
  \end{aligned}
  \end{equation}
The following can be seen by applying expanding \(\scv_n(\alpha,m)\) using Proposition~\ref{prop-combined-fna}~\ref{prop-combined-fna-l2} and noticing that the \(\frac{1}{m}\|\varphi_{\alpha}\|^2\) terms cancel.
\[
  \scv_n(\alpha,m) - R(\alpha,m) = T_1(\alpha) - T_2(\alpha),
\]
with
\[
  T_1(\alpha) = \frac{m-1}{nm(n-1)} \sum_{j=1}\sum_{k \neq j} \varphi_{\alpha}*\varphi_{\alpha}(\alpha)(Y_j-Y_k) - \frac2n \sum_{j=1}^n\mathcal{D}^{-1}f_{n,-j}^{\alpha}(Y_j)
\]
and
\[
  T_2(\alpha) = \frac{n+1}{n}\|\varphi_{\alpha}*g*f\|^2 - 2 \int [\varphi_{\alpha}*g*f(x)]f(x)\,dx.
\]

Applying this to the paired terms in Equation~\ref{eq-addsub} and then applying the triangle inequality, we find the upper bound
\begin{equation}
\label{eq-threeterms}
\begin{aligned}
  |\scv_n(\alpha',m) - R(\alpha',m)| &\leq |\scv_n(\alpha_j,m) - R(\alpha_j,m)|\\
  &\qquad + |T_1(\alpha') - T_1(\alpha_j)|\\
  &\qquad + |T_2(\alpha') - T_2(\alpha_j)|.
\end{aligned}
\end{equation}
Since \(\alpha_j \in \mathcal{A}_n\), by Lemma~\ref{lemma-finitesupdiff} there is a sequence \(C_{0,n} \to 0\) s.t. \(|\scv_n(\alpha_j,m) - R(\alpha_j,m)| \leq C_{0,n}r_m\), so now we only focus on the final two terms.

Now, since \(\alpha_j\) is the nearest element of \(\mathcal{A}_j\) to \(\alpha'\), we have \(|\alpha_j - \alpha'| \leq \frac{(\mu - \lambda)b_m}{k_n-1}\).
Using Lemma~\ref{lemma-T1T2bounds} to bound \(\frac{\partial}{\partial \alpha} T_2(\alpha)\),
\[\begin{aligned}
  |T_2(\alpha') - T_2(\alpha_j)| &= \left| \int_{\alpha'}^{\alpha_j} \frac{\partial}{\partial \alpha}T_2(\alpha)\,d\alpha \right|\\
                                 &\leq \int_{\alpha'}^{\alpha_j} \left| \frac{\partial}{\partial \alpha}T_2(\alpha)\right| \,d\alpha \\
                                 &\leq |\alpha_j - \alpha'| \sup_{\lambda b_m \leq \alpha \leq \mu b_m} \left| \frac{\partial}{\partial \alpha}T_2(\alpha)\right|\\
                                 &\leq \frac{(\mu - \lambda)b_m}{k_n-1} \times \frac{C_2}{\alpha}\\
                                 &\leq \frac{(\mu - \lambda)b_m}{k_n-1} \times \frac{C_2}{\lambda b_m}\\
                                 &\leq \frac{C'_2}{k_n-1}.
\end{aligned}\]
Proceding similarly for the \(T_1(\alpha)\) term,
\[\begin{aligned}
  |T_1(\alpha') - T_1(\alpha_j)| &\leq \frac{(\mu - \lambda)b_m}{k_n-1} \times \frac{C_1}{\alpha^2}\\
                                &\leq \frac{(\mu - \lambda)b_m}{k_n-1} \times \frac{C_1}{\lambda^2 b_m^2}\\
                                &\leq \frac{C'_1}{b_m(k_n-1)}.
\end{aligned}\]

At this point we have
\[
 |\scv_n(\alpha',m) - R(\alpha',m)|  \leq \left[ C_{0,n} + \frac{C_1}{r_mb_m(k_n-1)} + \frac{C_2}{r_m(k_n - 1)} \right] r_m,
\]
and \(C_{0,n} \to 0\), so if the other two terms in brackets vanish as \(n \to \infty\), then we are finished.
By Assumption~\ref{assumption:power}, \(n(r_mb_m)^2 = n^\delta\) with \(\delta > 0\).
This implies that \(\frac{C_1}{n^p(r_mb_m)^2} \to 0\) for any \(p \geq 1\); thus \(\frac{C_1}{n^{\frac{1}{2}}r_mb_m} \to 0\).
Since we have chosen \(k_n > n\), \(\frac{C_1}{r_mb_m(k_n-1)} \to 0\), as needed.
Finally, \(r_m > r_n\) since \(m < n\), and since \(r_n\) is a rate of convergence in deconvolution, we have \(n^{-1} = o(r_n)\).
Thus \(\frac{C_2}{r_m(k_n - 1)} \leq \frac{C_2}{r_n(k_n - 1)} \leq \frac{C_2}{r_nn}\to 0\).

\end{proof}

\begin{proof}[Proof of Corollary~\ref{cor-minimasupdiff}.]
  By Theorem~\ref{theorem-contsupdiff_as}, \(\sup_{\alpha \in \mathcal{A}_n} |\scv_n(\alpha) - R(\alpha,m)| \leq C_n r_m\) with \(C_n \to 0\).
  Write
  \[\begin{aligned}
   R(\alpha_m',m) - R(\alpha_m^{*},m) &=  R(\alpha_m',m) - \scv_n(\alpha_m')\\
                                     &\qquad  + \scv_n(\alpha_m') - \scv_n(\alpha_m^{*})\\
                                     &\qquad  + \scv_n(\alpha_m^{*}) - R(\alpha_m^{*},m).
  \end{aligned}\]
The first and third pairs are addressed by Theorem~\ref{theorem-contsupdiff_as}, while the second pair is negative because \(\alpha_m'\) minimizes \(\scv_n(\cdot)\), so
\[
   R(\alpha_m',m) - R(\alpha_m^{*},m) \leq  C_nr_m + 0 + C_n r_m = 2C_n r_m,
\]
i.e. \(R(\alpha_m',m) - R(\alpha_m^{*},m) = o(r_m)\).
\end{proof}

\begin{proof}[Proof of Corollary~\ref{cor-f2normest}]
  By Definition~\ref{def-r}, \(R(\alpha_m^{*},m) = \mathbb{E}\|f_n^{\alpha_m^{*}} - f\|^2 - \|f\|^2\), i.e. \(-R(\alpha_m^{*},m) = \|f\|^2 - \mathbb{E}\|f_m^{\alpha_m^{*}} - f\|^2\).
  Since \(\alpha_m^{*}\) minimizes \(\mathbb{E}\|f_m^{\alpha} - f\|^2\) in \(\alpha\), \(\mathbb{E}\|f_m^{\alpha_m^{*}} - f\|^2 = O(r_m)\).
  By Corollary~\ref{cor-minimasupdiff}, \(-\scv_n(\alpha_m',m) = -R(\alpha_m^{*},m) + o(r_m) = \|f\|^2 + O(r_m)\) almost surely, as needed.
\end{proof}

\subsection{Lemmas}\label{sec:appendix:lemmas}

\begin{lemma}\label{lemma-Qn}
  Let \(\tilde Q_n(t) = \frac{1}{n(n-1)}\sum_{j=1}^n\sum_{k \neq j} e^{-it(Y_j-Y_k)}\). Then \(\tilde Q_n(t)\) is an unbiased estimator of \(|\tilde h(t)|^2\), i.e. \(\mathbb{E}\tilde Q_n(t) = |\tilde h(t)|^2\).
\end{lemma}
\begin{proof}[Proof of Lemma~\ref{lemma-Qn}.]
  Since \(j \neq k\) in every summand and the \(Y_i\) are independent,
  \[\begin{aligned}
   \mathbb{E}\tilde Q_n(t)  &= \frac{1}{n(n-1)}\sum_{j=1}^n\sum_{k \neq j} \mathbb{E}e^{-it(Y_j-Y_k)}\\
                            &= \frac{1}{n(n-1)}\sum_{j=1}^n\sum_{k \neq j} \mathbb{E}[e^{-itY_j}]\mathbb{E}[e^{itY_k}]\\
                            &= \frac{1}{n(n-1)}\sum_{j=1}^n\sum_{k \neq j} \tilde h(t) \overline{\tilde h(t)}\\
                            &= |\tilde h(t)|^2
  \end{aligned}\]
\end{proof}

\begin{lemma}\label{lemma-multilinear}
  Let \(V\) be a vector space over the reals, and \(\mathcal{F} \subset V\) be convex with \(\vecspan(\mathcal{F}) \supset V\).
  Suppose \(M: V^n \to \mathbb{R}\) is multilinear and symmetric in its arguments.
  If \(M(u,u,\dots,u) = 0\) for all \(u \in \mathcal{F}\), then \(M(u_1,u_2, \dots u_n) = 0\) for all \(u_1, \dots, u_n \in V\).
\end{lemma}
\begin{proof}[Proof of Lemma~\ref{lemma-multilinear}]
  Since \(\vecspan(\mathcal{F}) \supset V\), it suffices to show that if \(u_1, \dots, u_n \in \mathcal{F}\), then \(M(u_1, \dots, u_n) = 0\).
  For brevity, for any symmetric multilinear form \(L: V^p \to \mathbb{R}\) we will write \(L_k(u,v) = L(u, \dots,v, \dots)\) with \(u\) appearing exactly \(k\) times and \(v\) appearing exactly \(p-k\) times. We will also write \(\ell(u) = L(u, u, \dots, u)\).
  We prove this by induction on the number \(j\) of unique arguments \(u_1, \dots, u_n\).
  The base case of \(j = 1\) is exactly the hypothesis of the lemma.

  For the inductive step, assume that if there are no more than \(j\) unique arguments, then \(M(u_1, \dots, u_n) = 0\).
  Fix \(j-1\) elements \(v_1, \dots, v_{j-1} \in \mathcal{F}\) and define a new symmetric, multilinear form by \(M^j(u_1, \dots, u_{n-j+1}) = M(v_1, \dots, v_{j-1}, u_1, \dots, u_{n-j+1})\) (if \(j = 1\), then we fix zero of the arguments, so \(M^j = M\)).
  Note that by the inductive hypothesis, \(m^j(u) = 0\) for any \(u \in \mathcal{F}\), since \(m^j(u)\) is \(M(u_1, \dots, u_n)\) evaluated with \(j\) or fewer unique arguments -- specifically the arguments \(v_1, \dots, v_{j-1}, u\).
  Fix \(u,v \in \mathcal{F}\), and set \(w_t = tu + (1-t)v\).
  Since \(\mathcal{F}\) is convex, \(w_t \in \mathcal{F}\), so \(m^j(w_t) = 0\).
  By the linearity and symmetry of \(M^j\), we can expand \(m^j(w_t)\) to find
  \[\begin{aligned}
   0 &= m^j(w_t)\\ 
     &= t^{n-j+1}m^j(u) + (1-t)^{n-j+1}m^j(v) + \sum_{k=1}^{n-1} \binom{n}{k}t^k(1-t)^{n-j+1-k}M^j_k(u,v)\\
     &= \sum_{k=1}^{n-1} \binom{n}{k}t^k(1-t)^{n-j+1-k}M^j_k(u,v),\\
  \end{aligned}\]
  where we used that \(m^j(u) = m^j(v) = 0\).
  The functions \(t^k(1-t)^{n-j+1-k}\) are linearly independent and the expression vanishes for all \(t \in [0,1]\), so the coefficients \(\binom{n}{k}M^j_k(u,v)\) all vanish; in particular \(M^j_k(u,v) = 0\) for all \(u,v \in \mathcal{F}\), and since \(v_1, \dots, v_{j-1}\) were also arbitrary, \(M(u_1, \dots, u_n)\) vanishes whenever there are no more than \(j+1\) unique arguments, as needed.

\end{proof}

\begin{lemma}\label{lemma-finitesupdiff}
  For \(\lambda b_{m} \leq \alpha_i \leq \mu b_{m}\), we have the bound
  \begin{equation}
    \mathbb{P}(\sup_{1 \leq i \leq k_n} r_{m}^{-1}|\scv_n(\alpha_i,m) - R(\alpha_i,m)| > \varepsilon) \leq 2k_n\exp\left \{-C\varepsilon^2n(r_{m}b_{m})^2\right \}.
  \end{equation}
  Thus if \((r_{m}b_{m})^{-2} = o(n)\) and \(k_n = o(e^{n(r_{m}b_{m})^2})\),
  \begin{equation}
    \sup_{1 \leq i \leq k_n} |\scv_n(\alpha_i,m) - R(\alpha_i,m)| = o_p(r_{m}).
  \end{equation}
\end{lemma}
\begin{proof}[Proof of Lemma~\ref{lemma-finitesupdiff}]
  We use McDiarmid's inequality.
  By the bound in Proposition~\ref{prop-ustat}, we have that \(|\theta_\alpha(x) - \theta_\alpha(y)| \leq C/\alpha\).
  Use the U-statistic representation of \(\scv_n(\alpha,m)\) to find, noticing that all terms not involving \(Y_i'\) vanish,
  \begin{equation}
    \begin{aligned}
    |\scv_n(\alpha,m;Y_1, &\dots, Y_i, \dots, Y_n) - \scv_n(\alpha,m;Y_1, \dots, Y_i', \dots, Y_n)| \\
    &\leq \left | \frac{2}{n(n-1)} \sum_{j \neq i} \theta_\alpha(Y_j - Y_i) - \theta_\alpha(Y_j - Y_i')\right |\\
    &\leq \frac{2}{n(n-1)} \sum_{j \neq i} |\theta_\alpha(Y_j - Y_i) - \theta_\alpha(Y_j - Y_i')|\\
    &\leq \frac{2}{n(n-1)} \cdot (n-1) \cdot \frac{C}{\alpha}\\
    &= \frac{C}{n\alpha}.
    \end{aligned}
  \end{equation}
  Thus we have the coordinatewise bounds \(c_i = \frac{C}{n\alpha}\) so that \(\sum_{i=1}^n c_i^2 = \frac{C^2}{n\alpha^2}\).
  Applying McDiarmid's inequality gives
  \[\begin{aligned}
    \mathbb{P}(r_{m}^{-1} |\scv_n(\alpha,m) - R(\alpha,m)| > \varepsilon) &\leq 2\exp\{-2C_1\varepsilon^2r_{m}^2n\alpha^2\}\\
                                                                            &\leq 2\exp\{-2C_2\varepsilon^2(r_{m}b_{m})^2n\}
  \end{aligned}\]
  A union bound gives the stated result.
\end{proof}

\begin{lemma}\label{lemma-T1T2bounds}
  With \(T_1(\alpha)\) and \(T_2(\alpha)\) defined as in the proof of Theorem~\ref{theorem-contsupdiff_as}, we have, for \(\alpha < M\),
  \[
    \frac{\partial}{\partial \alpha}T_1(\alpha) \leq \frac{C_1}{\alpha^2}\qquad\qquad\frac{\partial}{\partial \alpha}T_2(\alpha) \leq \frac{C_2}{\alpha}
  \]
  The constant \(C_1\) depends on \(M\), but the constant \(C_2\) does not.
\end{lemma}
\begin{proof}[Proof of Lemma~\ref{lemma-T1T2bounds}]
Taking Fourier transforms we have
\[\begin{aligned} %
 T_1(\alpha) &= \frac{1}{2\pi}\int \frac{m-1}{m}|\tilde \varphi_{\alpha}(t)|^2 \tilde Q_n(t)\,dt - \frac{1}{\pi} \int \frac{1}{|\tilde g(t)|^2 + \alpha |t|^{2\nu}}\tilde Q_n(t)\,dt\\
             &= \frac{1}{2\pi}\int \frac{\frac{m-1}{m}|\tilde g(t)|^2}{(|\tilde g(t)|^2 + \alpha |t|^{2\nu})^2} \tilde Q_n(t)\,dt - \frac{1}{\pi} \int \frac{1}{|\tilde g(t)|^2 + \alpha |t|^{2\nu}}\tilde Q_n(t)\,dt
\end{aligned}\]
Differentiating under the integral sign,
\[\begin{aligned}
 \frac{\partial}{\partial \alpha}T_1(\alpha) &= -\frac{1}{\pi}\int \frac{\frac{m-1}{m}|t|^{2\nu}|\tilde g(t)|^2}{(|\tilde g(t)|^2 + \alpha |t|^{2\nu})^3} \tilde Q_n(t)\,dt + \frac{1}{\pi} \int \frac{|t|^{2\nu}}{(|\tilde g(t)|^2 + \alpha |t|^{2\nu})^2}\tilde Q_n(t)\,dt\\
                                             &= \frac{1}{\pi}\int \frac{\alpha |t|^{4\nu} + \frac{1}{m}|t|^{2\nu}|\tilde g(t)|^2}{(|\tilde g(t)|^2 + \alpha |t|^{2\nu})^3} \tilde Q_n(t)\,dt.
\end{aligned}\]
Since \(|\tilde Q_n(t)| \leq 1\), we can bound this by
\begin{equation}\label{eq-T1bound}
  \begin{aligned}
  \frac{\partial}{\partial \alpha}T_1(\alpha) &\leq \frac{1}{\pi}\int \frac{\alpha |t|^{4\nu} + \frac{1}{m}|t|^{2\nu}|\tilde g(t)|^2}{(|\tilde g(t)|^2 + \alpha |t|^{2\nu})^3}\,dt\\
                                              &= \frac{1}{\pi}\int \frac{\alpha |t|^{4\nu}}{(|\tilde g(t)|^2 + \alpha |t|^{2\nu})^3}\,dt + \frac{1}{m\pi}\int \frac{|t|^{2\nu}|\tilde g(t)|^2}{(|\tilde g(t)|^2 + \alpha |t|^{2\nu})^3}\,dt\\
  \end{aligned}
\end{equation}
We'll show that \(\frac{\partial}{\partial \alpha}T_1(\alpha) \leq \frac{C_1}{\alpha^2}\).

For the first integral in the upper bound of Equation~\ref{eq-T1bound}, dropping the first and second term in the denominator of the integrand in turn yields upper bounds \(\alpha^{-2}|t|^{-2\nu}\) and \(\alpha|t|^{4\nu}|\tilde g(t)|^{-6}\) respectively.
Now choose \(\varepsilon > 0\) small enough that \(|\tilde g(t)| \geq \delta > 0\) for all \(0 \leq t \leq \varepsilon\).
Then integrate the first bound away from zero and the second bound near zero:
\[\begin{aligned}
\frac{1}{\pi}\int \frac{\alpha |t|^{4\nu}}{(|\tilde g(t)|^2 + \alpha |t|^{2\nu})^3}\,dt &\leq \frac{1}{\pi}\left [\alpha \int_{-\varepsilon}^{\varepsilon} t^{4\nu}|\tilde g(t)|^{-6}\,dt + \alpha^{-2}\int_{|t| > \varepsilon} |t|^{-2\nu}\,dt \right ]\\
                                                                                     &\leq \frac{1}{\pi} \left [ 2\alpha\varepsilon^{4\nu + 1}\delta^{-6} + \frac{2}{\alpha^2(2\nu-1)}\varepsilon^{1-2\nu} \right]\\
                                                                                     &\leq \frac{C_1'}{\alpha^2}.
\end{aligned}\]

For the second integral in the upper bound of Equation~\ref{eq-T1bound}, we take a similar approach.
The integrand is upper-bounded by \(\frac{4}{27\alpha^2}|t|^{-2\nu}\), which can be seen by applying the AM-GM inequality to \(|\tilde g(t)|^2 + \frac{\alpha}{2}|t|^{2\nu} + \frac{\alpha}{2}|t|^{2\nu}\) in the denominator.
On the other hand, dropping the \(\alpha |t|^{2\nu}\) in the denominator yields a different upper bound \(t^{2\nu}|\tilde g(t)|^{-4}\).
Now integrate these bounds away from and near zero, respectively:
\[\begin{aligned}
  \frac{1}{m\pi}\int \frac{|t|^{2\nu}|\tilde g(t)|^2}{(|\tilde g(t)|^2 + \alpha |t|^{2\nu})^3}\,dt  &\leq \frac{1}{m\pi} \left[ \int_{-\varepsilon}^{\varepsilon} t^{2\nu}|\tilde g(t)|^{-4} \,dt + \frac{4}{27\alpha^2}\int_{|t| > \varepsilon} |t|^{-2\nu}\,dt \right]\\
                                                                                                 &\leq \frac{1}{m\pi} \left[ 2\varepsilon^{2\nu+1}\delta^{-4} \,dt + \frac{8}{27\alpha^2(2\nu-1)}\varepsilon^{1-2\nu} \right ]\\
                                                                                                 &\leq \frac{C''_1}{\alpha^2}.
\end{aligned}\]
Combining these two, we have that \(\frac{\partial}{\partial \alpha}T_1(\alpha) \leq \frac{C_1}{\alpha^2}\).

Now, for \(T_2(\alpha)\), taking Fourier transforms,
\[
  T_2(\alpha) = \frac{1}{2\pi} \left[ \int \frac{\frac{n+1}{n}|\tilde g(t)|^4}{(|\tilde g(t)|^2 + \alpha |t|^{2\nu})^2} |\tilde f(t)|^2\,dt - 2 \int \frac{|\tilde g(t)|^2}{|\tilde g(t)|^2 + \alpha |t|^{2\nu}}|\tilde f(t)|^2\,dt \right]
\]
Differentiating under the integral,
\[
  \begin{aligned}
  \frac{\partial}{\partial \alpha}T_2(\alpha) &= \frac{1}{\pi} \int \frac{\alpha |t|^{4\nu}|\tilde g(t)|^2 - \frac1n |t|^{2\nu}|\tilde g(t)|^4}{(|\tilde g(t)|^2 + \alpha |t|^{2\nu})^3} |\tilde f(t)|^2\,dt\\
                                              &\leq \frac{1}{\pi} \int \frac{\alpha |t|^{4\nu}|\tilde g(t)|^2}{(|\tilde g(t)|^2 + \alpha |t|^{2\nu})^3} |\tilde f(t)|^2\,dt\\
  \end{aligned}
\]
It can be seen that \(\frac{\alpha |t|^{4\nu}|\tilde g(t)|^2}{(|\tilde g(t)|^2 + \alpha |t|^{2\nu})^3} \leq \frac{4}{27\pi \alpha}\) (apply the AM-GM inequality to \(|\tilde g(t)|^2 + \frac{\alpha}{2} |t|^{2\nu} + \frac{\alpha}{2} |t|^{2\nu}\)), so
\begin{equation}
\label{eq-T2bound}
\left[ \frac{\partial}{\partial \alpha}T_2(\alpha) \right] \leq \frac{8\|f\|^2}{27\pi}\alpha^{-1}.
\end{equation}

\end{proof}

\vskip .65cm
\noindent
David Kent, Cornell University
\vskip 2pt
\noindent
E-mail: dk657@cornell.edu
\vskip 2pt

\printbibliography

\end{document}